# Spéculation sur la géométrie en Égypte antique

©François Poisson, Québec, Canada

Résumé : Les grandes pyramides d'Égypte dissimulent des informations mathématiques ignorées jusqu'à aujourd'hui. Les mesures des trois grandes pyramides d'Égypte à Gizeh révèlent que les égyptiens de la IV$^e$ dynastie savaient calculer la circonférence, le volume et l'aire de la sphère 2000 ans avant les Grecs. La pyramide rhomboïdale affiche les trois grands problèmes de géométrie de l'antiquité : cubature de la sphère, duplication du cube, trisection de l'angle. Ainsi, nous affirmons que les Égyptiens voulaient graver de façon indestructible les concepts de base de la géométrie.

Abstract: The Great Pyramids of Egypt hide mathematic information unknown up to date. The measurements of the three Great Pyramids of Egypt at Giza show that Egyptians knew how to calculate the circumference, the volume and the area of the sphere 2000 years before Greeks. The Bent pyramid shows the three great problems of geometry of antiquity: cubature of the sphere, duplication of the cube, trisection of the angle. According to these findings, we assert that Egyptians wanted to engrave basic concepts in measures and positions of Great Pyramids.

## 1- Introduction

La majorité des égyptologues croient que les pyramides d'Égypte servent de tombeaux aux pharaons. Depuis des millénaires, ces immenses constructions soulèvent des mystères qui donnent lieu à d'innombrables théories et débats, autant amateurs que spécialistes. L'aspect mathématique et principalement géométrique fera l'objet des observations inédites et audacieuses qui vont suivre dans ce texte.

Plusieurs auteurs arabes affirmèrent que les Grandes Pyramides servaient de support pour enregistrer les connaissances de l'époque (Pochan[8], 1971, p.78). Rappelons que les grandes énigmes mathématiques de l'antiquité se rapportent aux caractéristiques des sphères, à doubler le volume d'un cube, à séparer un angle en trois parties égales, et à trouver des triangles rectangles dont les côtés sont formés de nombres entiers. De plus, les outils de l'époque requièrent de n'utiliser que la règle et le compas.

Les constructeurs des pyramides auraient-ils connu ces problèmes avant les Grecs de l'antiquité, voire même avant le papyrus de Rhind (Ahmes)? Pourquoi ne les retrouve-t-on pas dans l'architecture des monuments? Dans son livre, Corinna Rossi[9] (2003, p.68) rassemble une foule d'informations sur l'architecture des pyramides pour établir certains liens, traitant plusieurs aspects mathématiques. Elle conclut qu'il faut éviter d'utiliser des connaissances anachroniques pour interpréter les structures des Égyptiens et surtout qu'il faut détecter l'intention des constructeurs d'y inclure une connaissance spécifique.

Saccagées par le temps, les pyramides ont perdu leurs arêtes servant de points de référence aux mesures de toutes sortes. De plus, les outils de l'époque du Pharaon Kheops et les pierres immenses à déplacer obligeaient les constructeurs à certains compromis. En conséquence, les mesures de longueur et de distance représentent un grand défi pour reconstituer ces plans. Si on veut y lire un message, il devient très difficile de soutenir sa cause tant sur les écarts de mesures que sur les interprétations.

Nous désirons démontrer que les Égyptiens du temps de Kheops ont gravé les grandes énigmes de l'antiquité dans l'architecture des pyramides, ce qui va bien au-delà des traces écrites que nous possédons sur leurs connaissances.

| **Table 1. Gizeh** (Lehner) | **mètres** | **coudées** |
|---|---|---|
| Kheops base | 230,3 m | 440 c |
| Kheops hauteur | 146,6 m | 280 c |
| Khephren hauteur | 143,5 m | 274 c |
| Mykérinos hauteur | 65,0 m | 124 c |
| Mykérinos base NS | 104,6 m | 200 c |
| Mykérinos base EO | 102,2 m | 196 c |
| Djedefre Base | 104,6 m | 200 c |
| Djedefre Hauteur | 65,0 m | 124 c |

Table 1. Les mesures de Mark Lehner.

La coudée (c) représente l'unité de mesure de longueur à l'époque de la IV$^{e}$ dynastie des pharaons. Sa longueur équivaut à 0,5235 mètre (Gillings[2], 1982, p.220; Lehner[5], 1997, p108, Kheops : 230,33 m = 440 c). À moins d'avis contraire, nous utiliserons les mesures de Mark Lehner[5] (1997) dans ce texte telles que détaillées dans les tables 1 et 2.

## 1.1 Les racines carrées

John Legon[4] (1979) a observé qu'en prolongeant les côtés des pyramides de Kheops et de Mykérinos (figure 1), les longueurs des côtés perpendiculaires du grand rectangle correspondent à 1000 fois √2 et 1000 fois √3 coudées (1414 et 1732 coudées); ainsi la diagonale mesure 1000 fois √5 coudées (soit : 740,5 m, 907,0 m et 1170,8 m).

Gratuitement sur Internet, Google Earth[3] permet d'obtenir l'image du site de Gizeh. On peut tracer facilement le modèle théorique et prendre les mesures en mètres à l'aide du logiciel avec une précision de quelques mètres/km compte tenu de l'état des pyramides.

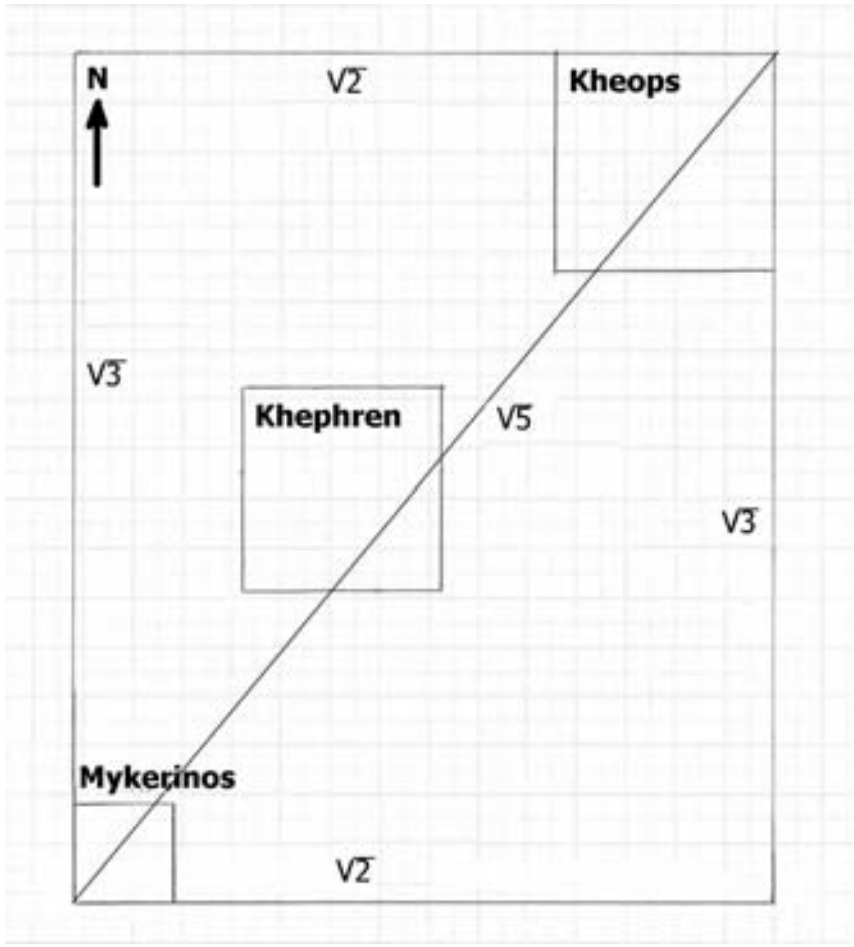

Figure 1. Site de Gizeh.

Rappelons que le système de fraction décimal ne s'est répandu qu'en 1585 par Simon Steven (Steichen[10], 1846) en publiant le DISME. Mais, un système semblable était utilisé au moyen âge pour représenter des racines carrées.

1.2 La circonférence.

Le mathématicien Paul Montel (Pochan[8], p.200) a reconnu la valeur de PI dans la pyramide de Kheops en trouvant que le périmètre de la base de Kheops équivalait à la circonférence du cercle dont le rayon serait la hauteur de la pyramide ($4C = 2\pi H$ ou 921,2 m = 921,1 m). Évidemment, il fut largement contesté prétextant le hasard et les chiffres sacrés, bien que la différence ne s'écarte que de 1/10 000 du périmètre.

De nombreux ouvrages font appel aux mathématiques pour donner toutes sortes d'interprétations aux mesures des pyramides. Nous y ajoutons cet article en ne traitant que l'aspect géométrique dans le but de démontrer que les Égyptiens possédaient ce savoir qui s'éteignit pendant 2000 ans pour des raisons inconnues.

## 2- Gizeh et la sphère

Non seulement les Égyptiens pouvaient reproduire la circonférence de la sphère, mais ils en exposent virtuellement le volume et la surface de la sphère à tous leurs visiteurs depuis des millénaires. Sur la figure 2, considérons la base de Kheops comme unité de longueur, de volume et de surface. Les hauteurs des pyramides de Gizeh servent de rayons à des sphères virtuelles : la circonférence de la sphère de Kheops, le volume de celle de Khephren et la surface de celle de Mykérinos concordent avec la base de Kheops.

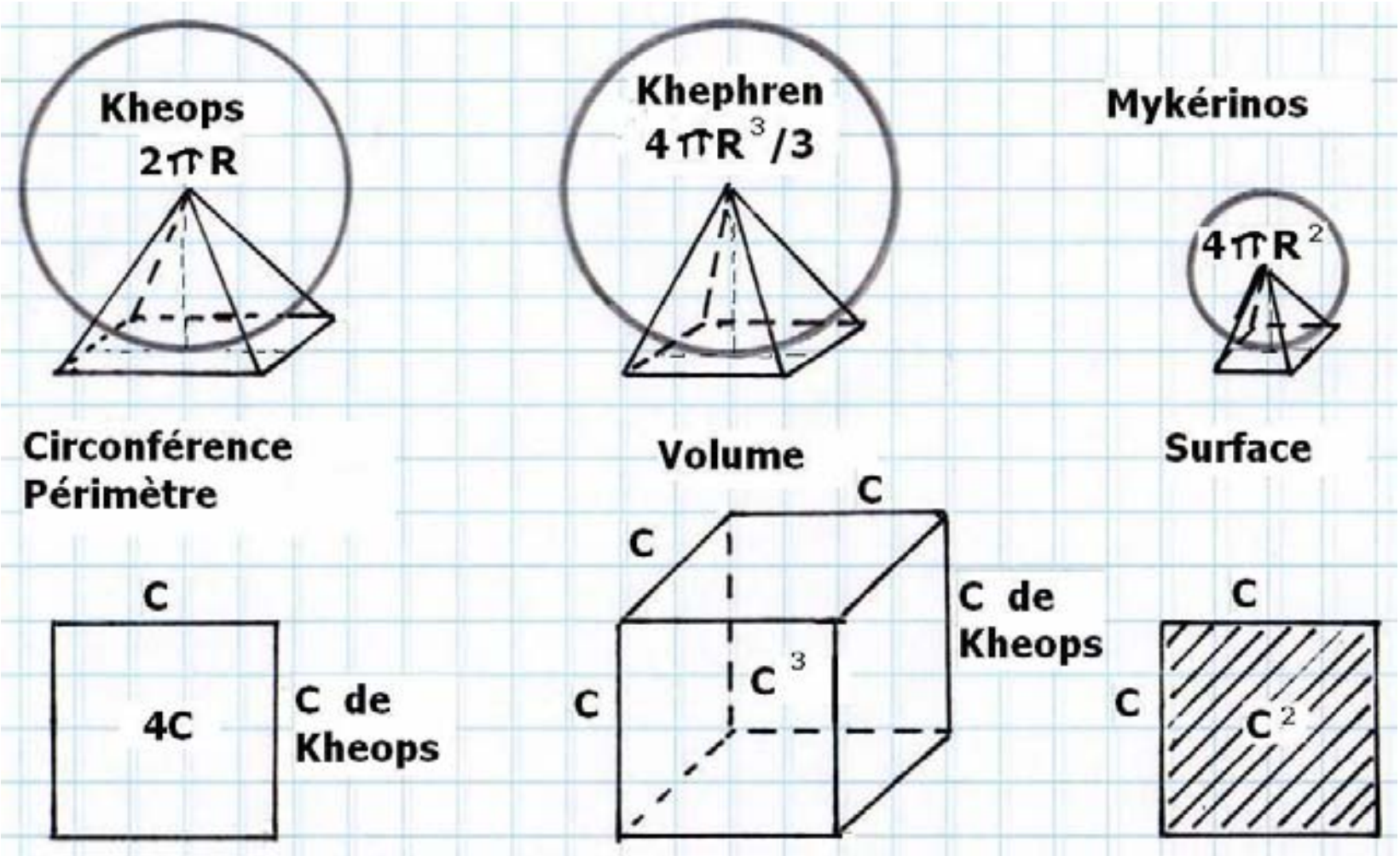

Figure 2. La hauteur de chaque pyramide correspond au rayon d'une la sphère.

2.1 Volume de la sphère de Khephren

Le volume du cube utilisant le côté de la base de Kheops est équivalent au volume de la sphère dont le rayon correspond à la hauteur de la pyramide de Khephren.

Le volume de la sphère $4\pi R^3/3$ mesure 12 377 824 $m^3$, soit à 1,3% du cube si on utilise la hauteur de 143,5 m (274 coudées), mais à une différence de volume de 5/10 000 si on utilise 142,9 m (273 c).

Nous suggérons d'utiliser 273 coudées pour la hauteur de Khephren mesurée par Maragioglio (Rossi[9], 2004, p.245). Précisons que Edwards (p.133) évalue la longueur des côtés de Khephren à 215,7 m (412 c) et que Lehner (p.122) l'évalue à 215 m (410,4 c) modifiant ainsi le résultat du calcul trigonométrique de la hauteur. Nous croyons que les Égyptiens voulaient obtenir une hauteur de 273 coudées en utilisant un côté de 410 coudées (214,7 m) et une pente de 4/3 (53°08').

2.2 L'aire de la sphère de Mykérinos

L'aire du carré utilisant le côté de la base de Kheops équivaut à l'aire de la sphère dont le rayon correspond à la hauteur de la pyramide de Mykérinos (65 m, 124 c).

L'aire de la sphère $4\pi R^2$ mesure 53 092 $m^2$, soit une différence de surface de 1/1000 relativement au carré de la base de Kheops (230,3 m).

En plus de comparer la surface de sa sphère avec la base de Kheops, Mykérinos établit un lien de volume et de circonférence avec sa propre base rectangulaire selon le côté utilisé.

- Le volume du cube utilisant le côté de la base de 104,6 m de Mykérinos est équivalent au volume de la sphère dont le rayon correspond à la hauteur 65 m de cette même pyramide, avec une différence de 5/1000 du volume.

- Le périmètre de la base de 102,2 m (x4) de Mykérinos équivaut à la circonférence du cercle dont le rayon correspond à la hauteur 65 m de cette même pyramide avec une différence de 1/1000 du périmètre.

2.3 Djedefre pourrait être magique aussi

Un autre pharaon de la quatrième dynastie a construit sa pyramide à Abouroash au nord-est de Gizeh. Ses proportions viennent d'être précisées en dessablant sa base : 106 m à la base et entre 57 et 67 m de hauteur (Nordland[6], 2008). Comme celle de Mykérinos, elle pourrait symboliser une sphère ayant un volume égal au cube de sa base si on utilisait 200 coudées (104,7 m) pour sa base et 124 coudées (65 m) de hauteur.

## 3- La pyramide rhomboïdale et la géométrie

La pyramide rhomboïdale fut construite par Snéfrou, premier pharaon de la quatrième dynastie et père de Kheops. Les égyptologues expliquent sa forme irrégulière par des corrections de la structure nécessitées par la fragilité du matériel utilisé.

Nous percevons plutôt une construction scientifique stupéfiante pour exprimer les trois grands problèmes de géométrie de l'antiquité à la vue de tous ses visiteurs :

- la duplication du cube : trouver la longueur de l'arête d'un cube qui aura le double du volume d'un autre cube;
- la cubature de la sphère : déterminer le rayon d'une sphère ayant le même volume qu'un cube déterminé;
- et la trisection de l'angle : séparer un angle en trois angles égaux.

La difficulté définie par les Grecs de l'antiquité à ces grands problèmes provient qu'il ne faut utiliser que la règle et le compas pour obtenir le résultat.

| **Table 2. Rhomboïdale** (Lehner) | **mètres** | **coudées** |
|---|---|---|
| Pyramide satellite base | 53 m | 101 c |
| Pyramide satellite hauteur | 32,5 m | 62 c |
| Pyramide rhomboïdale base | 188 m | 360 c |
| Pyramide rhomboïdale hauteur | 105 m | 200 c |

Table 2. Mesures de Lehner sur le site de la pyramide rhomboïdale.

3.1 Un cube dont le volume équivaut au double d'un autre cube.

La pyramide Rhomboïdale contient deux passages souterrains, suggérant qu'il s'agit de deux pyramides en une seule comme le présente la figure 3. Remarquons que la hauteur de la plus basse pyramide mesure 105 m (table 2) et pourrait représenter l'arête du cube de volume égal à 1 157 625 mètres cubes.

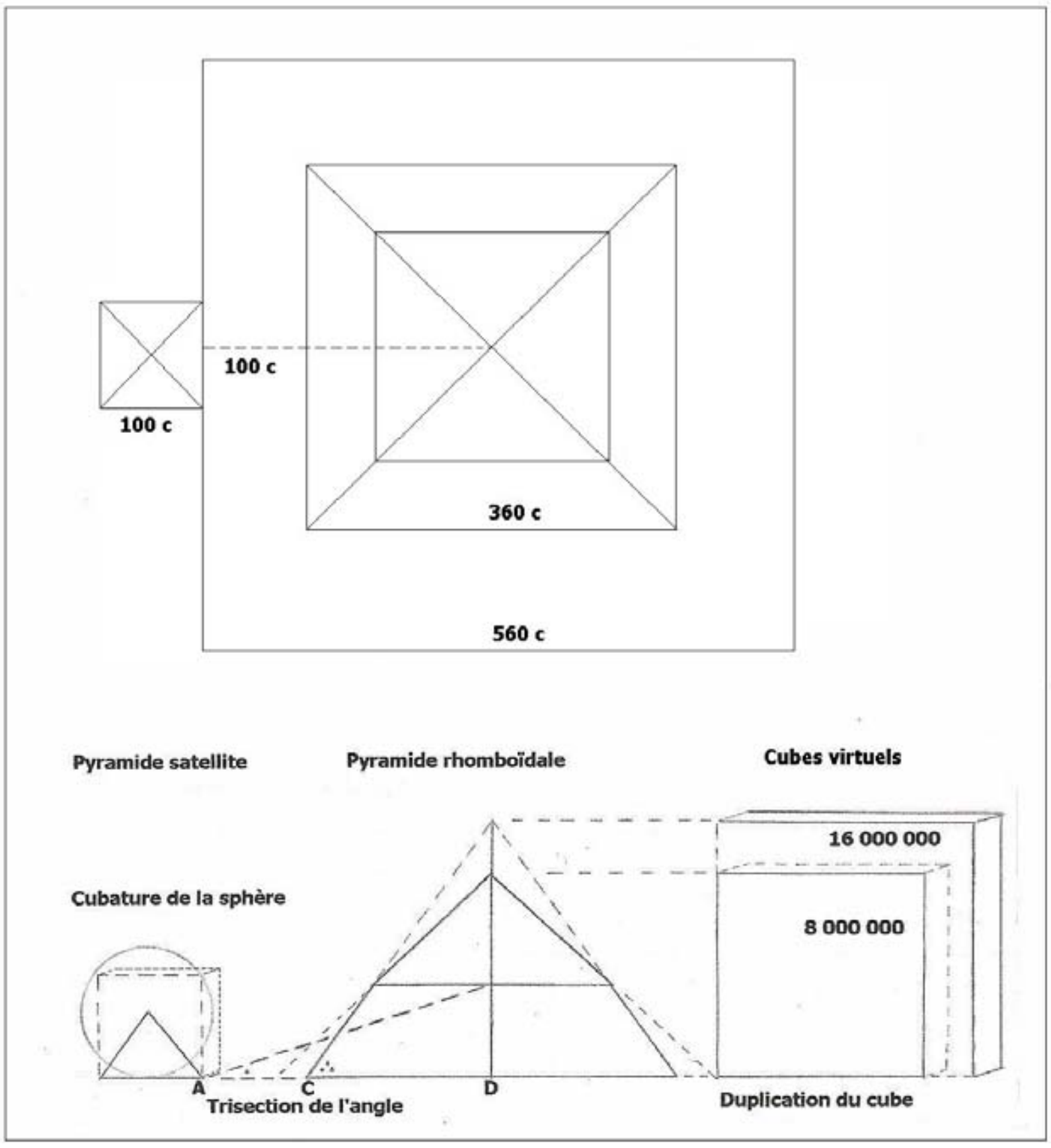

Figure 3**.** La pyramide rhomboïdale et les grands problèmes de l'antiquité**.**

La hauteur de la plus haute pyramide dont on prolonge les faces de ses côtés se calcule en utilisant la moitié de sa base (94 m, projections des côtés du triangle) multipliée par la pente 1,4 de son côté pour obtenir 131,6 m. Elle pourrait représenter l'arête du cube de 2 279 122 mètres cubes, soit le double du premier cube à 1,6% près.

Examinons le cube en utilisant la coudée égyptienne. La hauteur de la plus basse pyramide mesure 200 coudées (105 m) et représente l'arête du cube de 8 millions de coudées cubes. La hauteur de la plus haute pyramide dont les faces de ses côtés sont prolongées, mesure 252 coudées (131,9 m) et représente l'arête du cube de 16 millions de coudées cubes avec la précision de 2 sur 10 000. C'est la duplication du cube.

## 3.2 Technique égyptienne pour doubler le volume d'un cube**.**

Nous croyons que les Égyptiens suggèrent une solution à la duplication du cube par la position de ses monuments. Des indices permettent de trouver leur solution du problème. Considérons un nombre remarquable de 100 coudées (51.92 m, Pétrie[7], 1888, p.27) plutôt que 104,8 coudées (54,86 m, Edwards, p.84) pour le segment de droite AC entre la pyramide satellite et la pyramide rhomboïdale.

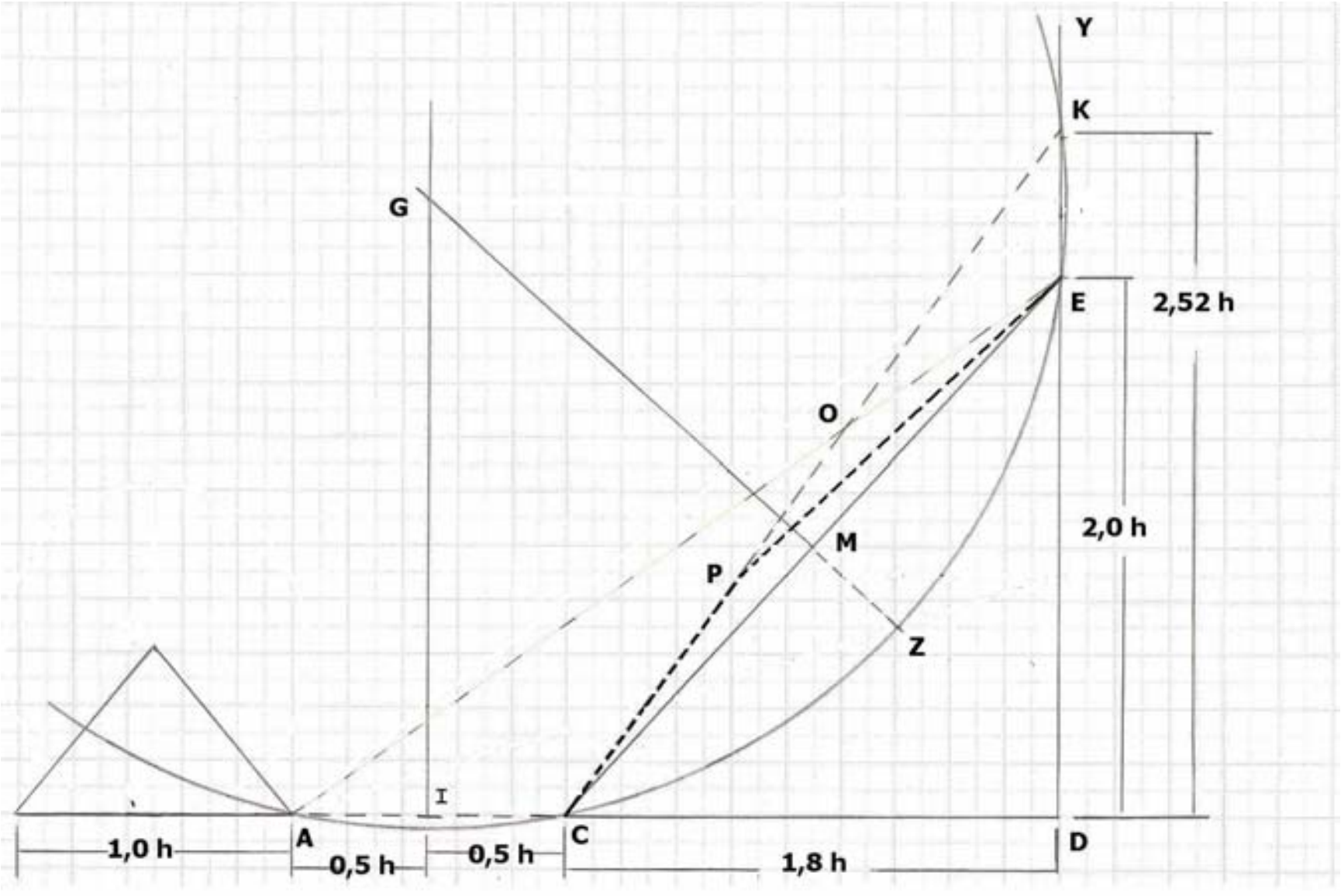

Figure 4. La solution des Égyptiens de la duplication du cube.

Pour trouver l'arête DK d'un cube de volume double d'un premier cube dont l'arête DE mesure 2h unités (figure 4) :

1- Élevons en D une perpendiculaire DY contenant une longueur DE de 2h unités.
2- CD représente les 1,8h unités de la demie base de la pyramide obtenue par des opérations simples de géométrie.
3- Situons AC la distance de 1h unités qui sépare la pyramide satellite de la grande pyramide dont le centre reçoit en I une perpendiculaire IG.
4- Joignons les points C et E et élevons sur CE une médiatrice MG où G est le point de rencontre avec la perpendiculaire IG.
5- Traçons le cercle de centre G passant par A, C, E et coupant DY en K.
6- Le volume $(DK)^3$ ou $(2{,}52h)^3$ sera très près du double du volume $(DE)^3$ ou $(2h)^3$ en unités cubes.

Voici les indices qui conduisent à cette proposition :

- La bande de 100 coudées entre la pyramide satellite et la grande pyramide permettant d'ériger une médiatrice au centre de l'unité semblable à la solution proposée par Nicodème, mais utilisant une relation simple plutôt qu'une moyenne prortionnelle.
- Les segments CK et AE qui se croisent en un point O ayant la même propriété que des cordes qui se croisent dans un cercle, soit CO x OK égale AO x OE.
- Les deux triangles semblables KCD et AED possèdent chacun un angle droit et l'angle CAE égale l'angle CKE car ils supportent le même arc CE. Les côtés sont proportionnels : DK/AD = CD/DE, ainsi DK/2,8h = 1,8h/2h, et enfin DK = 2,52h.

La disposition exceptionnelle des angles et des mesures ne semble pas être le fruit du hasard. Il s'agit d'une approximation ayant une précision de 2 sur 10 000. Ce qui pouvait être considéré comme suffisamment précis pour les Égyptiens.

3.3 Une sphère de même volume qu'un cube

Selon les mesures des égyptologues (table 2), la base de la pyramide satellite près de la pyramide rhomboïdale mesure 53 m. Considérons un cube dont l'arête de 53 m est située sur la même base que la pyramide satellite, son volume serait de 148 877 $m^3$. Imaginons une sphère dont le rayon correspond à la hauteur de la pyramide satellite de 32,5 m, son volume serait de 143 793 $m^3$. La différence de volume entre le cube et la sphère serait de 3,4%.

Passons en mode Égyptien. La base de la pyramide satellite mesure 53 m, ce qui donne 101 coudées en utilisant la mesure de la coudée (0,5235 m). Compte tenu de la détérioration de cette pyramide avec le temps, il y a de fortes chances pour que sa mesure fût à l'origine de 100 coudées (52,35 m, Pétrie[7], p.27). Si on conserve sa hauteur d'origine à 62 coudées (32,5 m), le volume de la sphère s'approche à 1,7/1000 de celui du cube. C'est la cubature de la sphère. Cette connaissance aurait été réutilisée quelques générations plus tard par les Égyptiens pour établir la hauteur de la pyramide de Khephren et celle de Mykérinos en fonction d'un autre volume.

3.4 Pour obtenir le rayon de la sphère

Sur la figure 5, utilisons le centre W de EK pour tracer le segment WA. Partant de WA, traçons la perpendiculaire AU qui rencontrera la médiatrice de la base de la pyramide satellite au point U, déterminant ainsi le rayon UV de la sphère (r = 50x280/226 = 62).

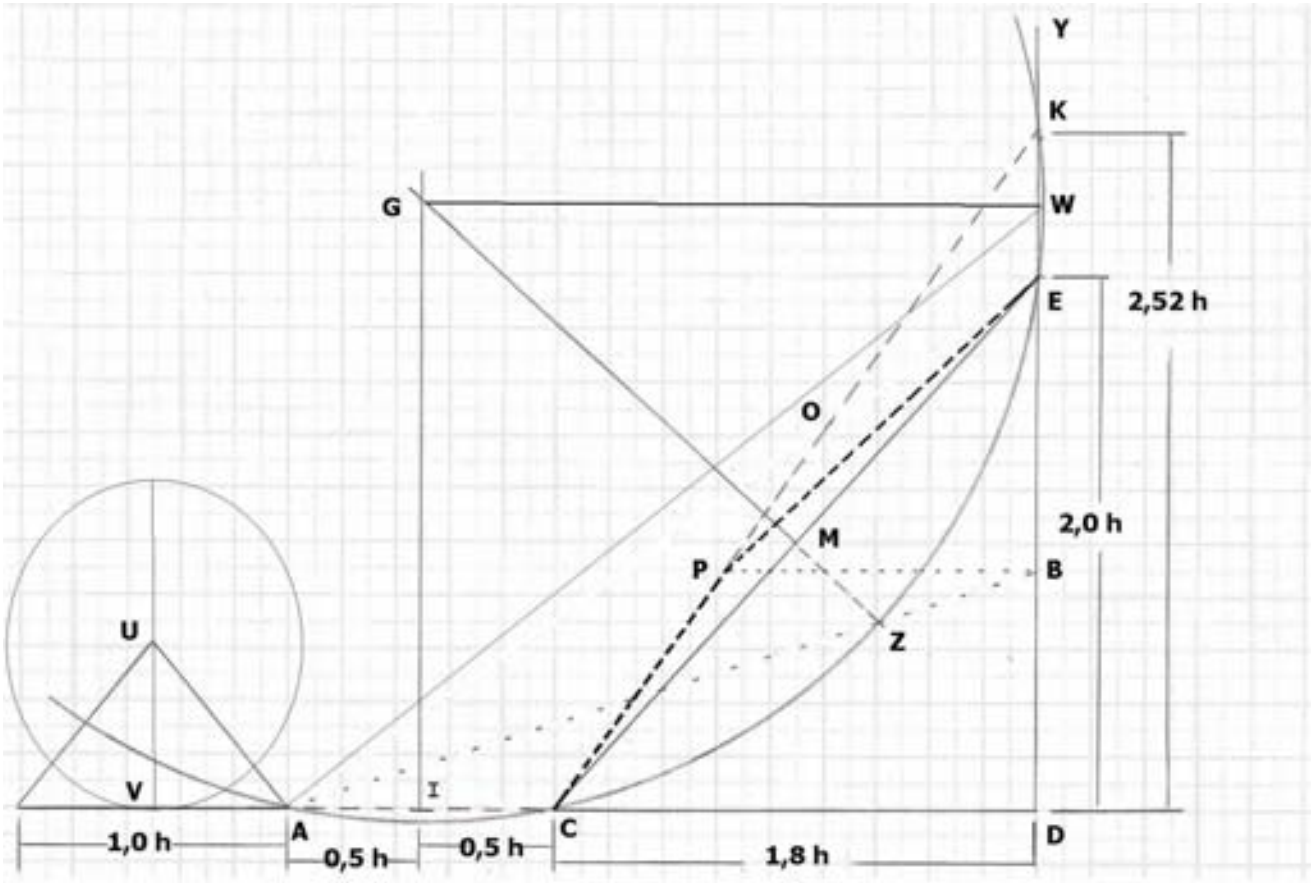

Figure 5. La solution des Égyptiens de la cubature de la sphère.

Nous ne pouvons pas expliquer la logique géométrique qui s'y rattache. Mais, cette coïncidence permet d'obtenir le rayon de la sphère de même volume que le cube de la base de cette pyramide satellite avec une précision de 1 sur 1000.

3.5 Un angle valant le tiers d'un autre angle

La face inférieure de la pyramide possède une pente de 1,4 mesurée avec une grande précision à 54°27'44'' par Lehner (p.102). La surface supérieure de la pyramide possède une pente de ayant un angle de 43°22'. Sur la figure 5, on peut retrouver l'angle stratégique KCD et son tiers BAD entre le segment AD longeant le sol et le segment AB reliant la base de la pyramide satellite et le centre du palier de la grande pyramide.

En utilisant 252 coudées comme longueur de DK et 180 comme longueur de CD, et selon les deux angles cités précédemment, on estime la hauteur du palier BP par calcul à 92,15 coudées (48,35 m). On obtient ainsi un angle BAD de 18°13'12'', comme étant le tiers du grand angle de la pyramide au lieu du tiers théoriquement calculé de 18°9'15'', qui situerait le palier à 91,81 coudées, soit une précision angulaire de 3,5 sur 1000. Ceci représente la trisection de l'angle.

S'il n'y avait que la duplication du cube à tenir compte, ce palier pourrait se situer à différentes hauteurs. Il suffirait de déplacer le point P le long de la ligne CK. Alors, la hauteur du palier a-t-elle été choisie stratégiquement pour que l'angle BAD mesure exactement le tiers de l'angle KCD ou s'agit-il d'un hasard?

Sur la figure 5, le point Z marque la rencontre du prolongement de GM avec le cercle et pourrait servir à localiser le point B représentant la hauteur du palier. Par trigonométrie, on calcule l'angle ZAD à 17,77°, soit de 2% inférieur, et le palier à 89,7 coudées. C'est loin de la précision que les Égyptiens ont démontrée jusqu'ici. Nous associons la trisection à la technique suivante pour séparer un angle en trois parties égales.

### 3.6 Technique égyptienne pour séparer un angle en trois parties égales

Encore ici, des indices sur le site de la pyramide rhomboïdale semblent nous suggérer une solution de la trisection d'un angle. Nous n'avons pas retrouvé ce procédé dans les documents de mathématiques consultés. Il s'agit d'une bonne approximation plus simple que la courbe conchoïde utilisée par les Grecs de l'antiquité.

Sur la figure 6, construisons un angle égal au tiers de l'angle KCD qui occupe une situation stratégique sur le site de la pyramide rhomboïdale.

1- Le segment CD représente la base de l'angle KCD.
2- Sur le segment AD, situons AC de une unité qui sépare la pyramide satellite de la grande pyramide. Le centre de AC reçoit en I une perpendiculaire IG.
3- Marquons S sur CK de sorte que CS = AI = CI = 0,5 unité.
4- Traçons SA qui croise IG en R.
5- Reportons la longueur CR = $CT_1$ sur CK produisant le point $T_1$ et traçons $T_1A$ qui coupe IG au point $X_1$, donnant une première estimation $T_1AI$ de l'angle recherché.
6- Si on veut plus de précision, reprenons l'étape 5 en reportant $CX_1$ sur CK produisant un nouveau $T_2$ que l'on relie au point A pour trouver un nouvel $X_2$, jusqu'à ce que $X_n$ et $T_n$ ne se déplacent que de façon infime.
7- L'angle TAI prolongé en BAD est le tiers de l'angle TCD ou KCD;

La configuration d'Archimède apparaît quand la limite de l'angle est atteinte :
XAI = XCI= a, car IX est sur une médiatrice;

TXC = 2a, car l'angle au centre d'un arc de cercle est le double d'un angle de même arc dont le sommet est situé sur la circonférence;
TCL = 2xTXC = 4a, même raison puisqu'on a deux cercles égaux;
DCL = a, opposé par le sommet à XCI;
TCD = TCL – DCL = 3a; par soustraction.

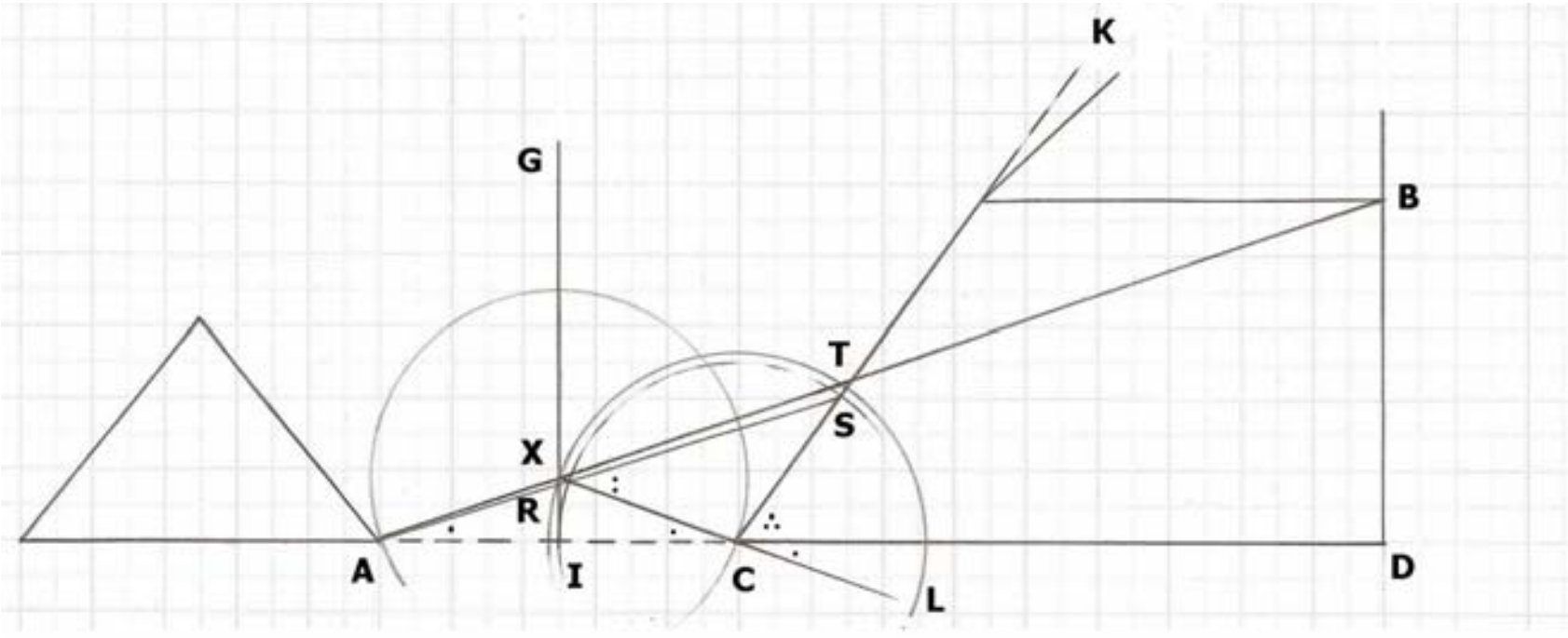

Figure 6. La solution des Égyptiens de la trisection de l'angle.

Cette solution utilise une répétition de la dernière opération pour obtenir le tiers de l'angle KCD de 54,46222°. En se référant à la table 3, la première estimation de $T_1AC$ donne un écart de 3 minutes, soit une erreur de 3/1000, alors que $T_2AC$ serait 13 fois plus précis (2/10 000) pour obtenir le tiers de l'angle calculé à 18,15407° .

| Table 3. Calcul trigonométrique de chaque itération de la trisection de l'angle. | | | | | | |
|---|---|---|---|---|---|---|
| Cycle | IXn | CXn | hauteur Tn | Proj.de ATn | Tiers angle | Cycle |
| | mètres | mètres | mètres | mètres | degrés | |
| S | 0,00000 | 50,00000 | 40,68662 | 129,06198 | 17,49739 | S |
| T1 | 15,76243 | 52,42570 | 42,66050 | 130,47190 | 18,10623 | T1 |
| T2 | 16,34854 | 52,60489 | 42,80631 | 130,57605 | 18,15052 | T2 |
| T3 | 16,39133 | 52,61821 | 42,81714 | 130,58379 | 18,15381 | T3 |
| T4 | 16,39451 | 52,61920 | 42,81795 | 130,58436 | 18,15405 | T4 |

Table 3. Calcul trigonométrique de chaque itération de la trisection de l'angle.

Le point commun de ces trois problèmes célèbres consiste à utiliser la perpendiculaire au centre de l'espace entre les deux pyramides. Cela parait très habile des Égyptiens de faire coïncider ce triple message pour démontrer la puissance de leur connaissance. De plus, si le pyramidion associé incorrectement à la pyramide rouge selon Rossi[9] (p. 207) était ramené au point stratégique I du site de la pyramide rhomboïdale, cela appuierait cette théorie : les trois problèmes célèbres rassemblés à un seul endroit.

## 4- Les pyramides et Pythagore

Le problème mathématique de la pyramide rouge, érigée également par Snéfrou, invite les pythagoriciens à trouver un triangle rectangle formé de nombres entiers et possédant des côtés perpendiculaires consécutifs. Le premier triangle étant très connu : $3^2+4^2=5^2$, quel est le suivant?

Énoncé par Alexandre Varille (Rossi[9], p.219), la réponse $20^2+21^2=29^2$ se trouve dans les dimensions de la pyramide rouge avec une amplification de dix, soit 200 coudées (104,7 m) de hauteur, 210 coudées (110 m) de base au sol et une hypoténuse de 290 coudées (151,8 m) (figure 7). La précision est de 3/1000 selon les mesures de Lehner (p.104).

Ce concept pourrait provenir de la construction de la pyramide rhomboïdale. En effet, les architectes auraient remarqué des multiples du triplet (20, 21, 29) en réalisant une pente presque identique à la dernière étape de la construction. Et cette trouvaille les impressionna suffisamment pour l'immortaliser dans la pyramide rouge.

Le troisième triangle de ce type possède des côtés mesurant $119^2+120^2=169^2$. On serait tenté de croire que la pyramide de Djoser représente ce triplet si on utilise les mesures de Edwards (p.35), soit 119 coudées (62,2 m) de hauteur et 120 coudées (62,6 m) de demie longueur du plus grand côté de sa base (figure 7). Mais, des mesures plus récentes de Lehner (p.85) fixent la hauteur à 114,4 coudées (60 m) et le demi grand côté à 115,5 coudées (60,5 m), ce qui demande une vérification avisée pour le confirmer ou le nier.

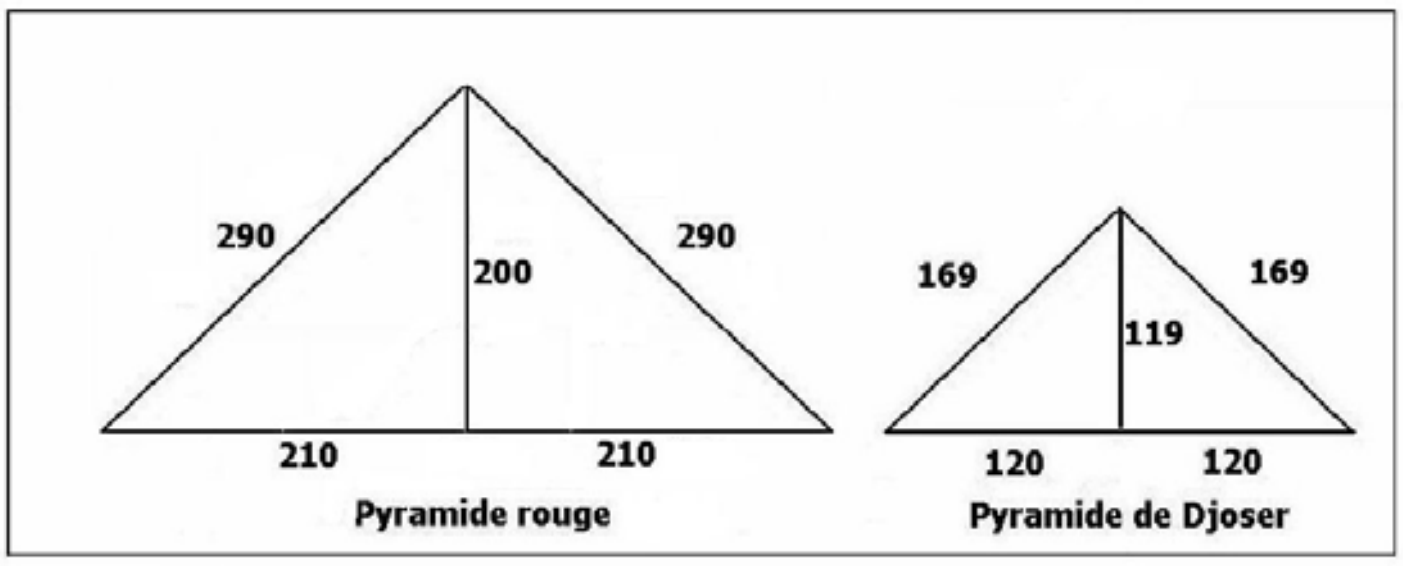

Figure 7. Les triangles pythagoriciens à côtés consécutifs.

Le mur entourant la pyramide de Djoser mesure 3142 coudées (1645 m) selon Lehner (p.84), soit très près de 1000 fois PI. Mais, s'agit-il d'un hasard?

La diagonale du rectangle de la base de la pyramide de Djoser mesure légèrement plus que 100 fois PI selon Edwards ($208^2+240^2=317,6^2$ coudées ou 166,2 m), mais légèrement moins selon Lehner ($208^2+231^2=310,8^2$ coudées ou 162,7 m). Ce qui demande une vérification avisée, car la moyenne donne 314,2 coudées, très près de 100 fois PI. Toutefois, la recherche de PI dans la diagonale semble en conflit avec la recherche du triplet (119, 120, 169), à moins de considérer de façon distincte le point de rencontre de la pente avec le sol et le côté de la pyramide. Cette partie révèle plus de spéculation.

Devant les formes irrégulières et détériorées de la pyramide de Djoser, il ne semble pas possible de lever un tel doute. Cela exigerait une grande précision impossible face à autant de détérioration du monument.

## 5- Conclusion

La connaissance de la géométrie à l'époque de Kheops soulève l'intérêt et les passions de nombreux amateurs et spécialistes des pyramides car aucun document de cette époque n'en fait état. Le site de Gizeh présente de toute évidence les caractéristiques des sphères.

La pyramide rhomboïdale représente de façon si détaillée et précise les trois grands problèmes de géométrie de l'antiquité que ce ne peut pas être le fruit du hasard. Ce site semble une œuvre d'art, telle une peinture ou une statue, exhibant virtuellement la cubature de la sphère, la duplication du cube et la trisection de l'angle. De plus, les Égyptiens auraient laissé des indices pour permettre de retrouver leurs solutions.

La pyramide rouge traite des nombres entiers de Pythagore. Plus loin dans le temps, la pyramide de Djoser montre sans clarté évidente la valeur de PI et le plus grand triangle connu à l'époque dont les côtés perpendiculaires sont des nombres entiers consécutifs.

Nous croyons que les Égyptiens ont laissé une marque chiffrée de leur savoir via les pyramides. Les lignes, les angles et les formes dénotent une intention claire de laisser un message pertinent dans leur architecture qui remet en question l'histoire des mathématiques au temps des pyramides.

## 6- Références